\documentclass[12pt,a4paper,twoside]{article}
\usepackage{amsmath,amsfonts,amssymb,a4,enumerate,epsfig,array,theorem}

\newfont{\eightrm}{cmr8}
\newfont{\fiverm}{cmr5}
\newfont{\bigrm}{cmr17}

\newcommand\tr{\operatorname{tr}}

\newcommand{\half}{{0.5}}

\newcommand{\ZZ}{{\mathbb{Z}}}

\newtheorem{prop}{Proposition}
\newtheorem{coro}[prop]{Corollary}
\newtheorem{theo}[prop]{Theorem}

\begin{document}

\title{Two results on tilings of quadriculated annuli}
\author{Nicolau C. Saldanha and Carlos Tomei}
\date{}
\maketitle


\begin{abstract}
We provide a more informal explanation of two results in our manuscript 
{\it Tilings of quadriculated annuli}.
Tilings of a quadriculated annulus $A$ are counted according to {\em volume}
(in the formal variable $q$) and {\em flux} (in $p$).
The generating function $\Phi_A(p,q)$ is such that,
for $q = -1$, the non-zero roots in $p$ are roots of unity
and for $q > 0$, real negative.
\end{abstract}

There is an unexpected rigidity for the $-1$-counting of tilings of
quadriculated disks, as described in \cite{DT}:

\begin{theo}
\label{theo:DT}
Let $D$ be a quadriculated disk.
Then the determinant of the adjacency matrix of the squares of $D$
equals -1, 0 or 1.
\end{theo}  

In \cite{ST} we prove Theorem \ref{theo:ST} below which, in a sense,
extends this rigidity to annuli.
In this shorter text, we present a more informal proof of this result.

\section{Connectivity}

We assume that the squares in a quadriculated annulus $A$
are colored in a checkerboard pattern
and that the numbers of black and white squares are equal.
Without loss, $A$ is embedded in the plane.
\footnote{1991 {\em Mathematics Subject Classification}.
Primary 05B45, Secondary 05A15, 05C50, 05E05.
{\em Keywords and phrases} Quadriculated surfaces, tilings by dominoes,
dimers.}
\footnote{The authors gratefully acknowledge the support of
CNPq, Faperj and Finep.}

\begin{figure}[ht]
\begin{center}
\epsfig{height=3cm,file=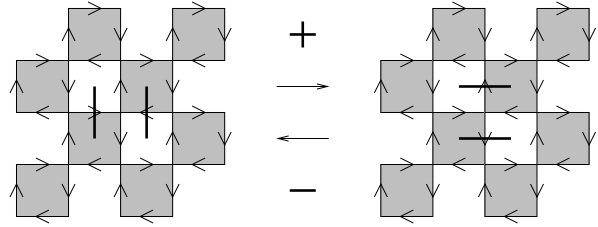}
\end{center}
\label{fig:signflip}
\caption{Positive and negative flips}
\end{figure}

Let $T_A$ be the set of domino tilings of $A$.
Two tilings are {\it adjacent} if they differ by a {\it flip},
a $90^\circ$ rotation of two dominoes filling a $2 \times 2$ square.
Assign the counterclockwise (resp. clockwise) orientation
to the boundary of white (resp. black) squares.
A flip is {\it positive} (resp. {\it negative})
if in the original $2 \times 2$ square
the two sides of squares from center to boundary not trespassing
dominoes point outwards (resp. inwards), 
as in Figure \ref{fig:signflip}.
The attribution of signs to flips is {\it exact}
in the sense any closed sequence of flips contains the same number
of positive and negative flips
(this follows easily from properties of height functions and sections
as discussed in \cite{STCR};
height functions were introduced by Thurston in \cite{T}).

\newpage

\goodbreak

\begin{figure}[ht]
\begin{center}
\epsfig{height=3cm,file=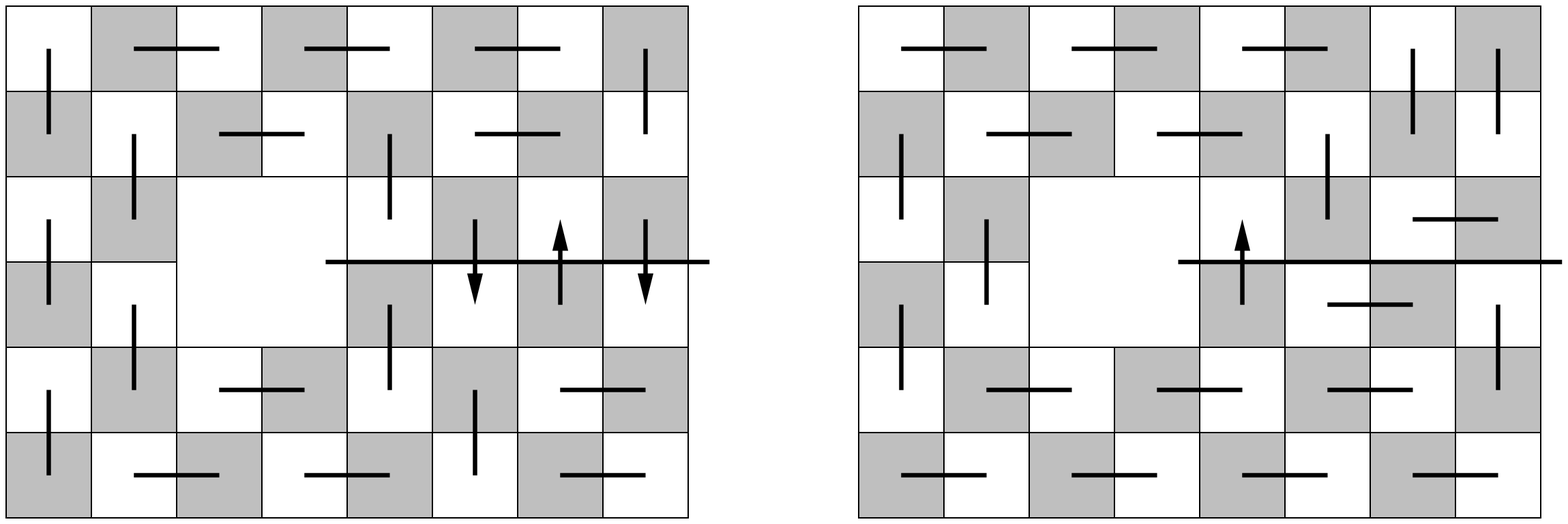}
\end{center}
\label{fig:miniflux}
\caption{Tilings with counterclockwise fluxes $-1$ and $1$}
\end{figure}

When are two tilings joined by a sequence of flips?
They clearly must have the same {\it flux} across a {\it cut}:
in Figure \ref{fig:miniflux}, the two tilings have different fluxes.
Indeed, flips don't alter flux.

\begin{figure}[ht]
\begin{center}
\epsfig{height=27mm,file=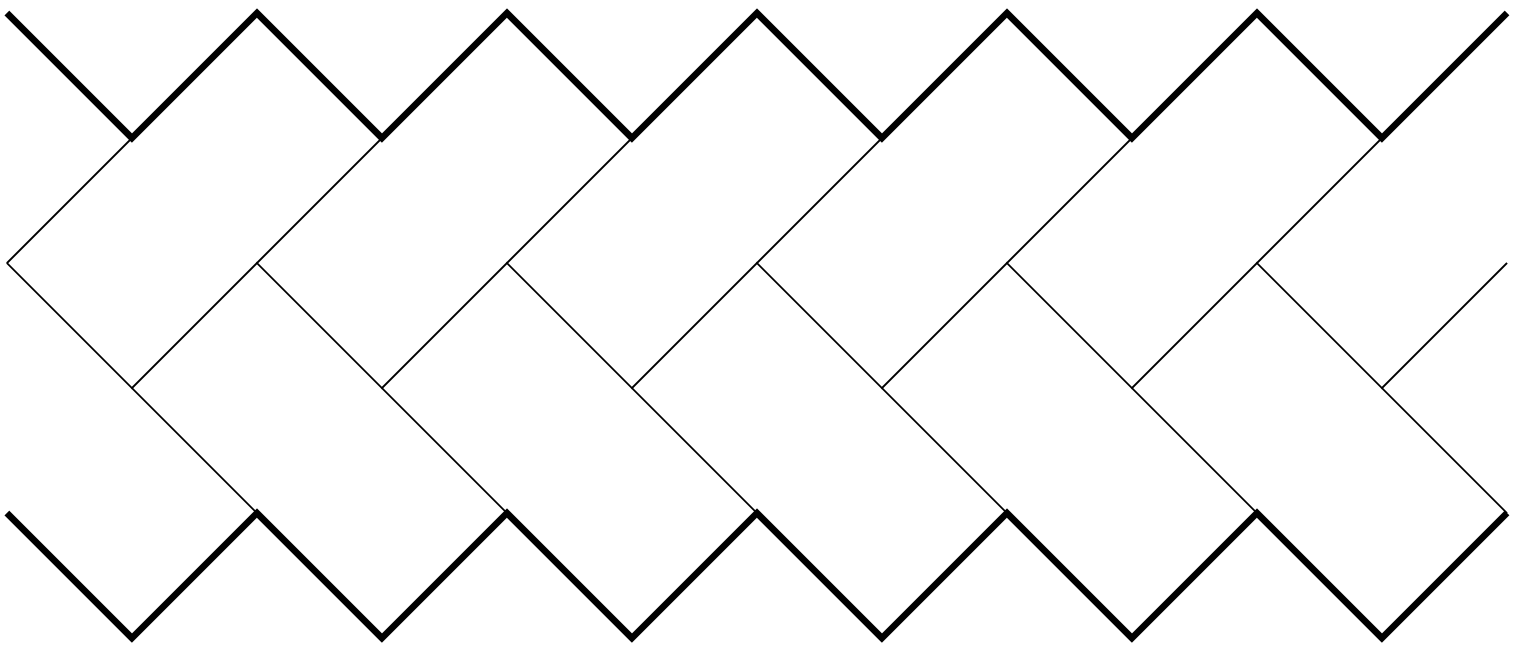}
\quad \epsfig{height=27mm,file=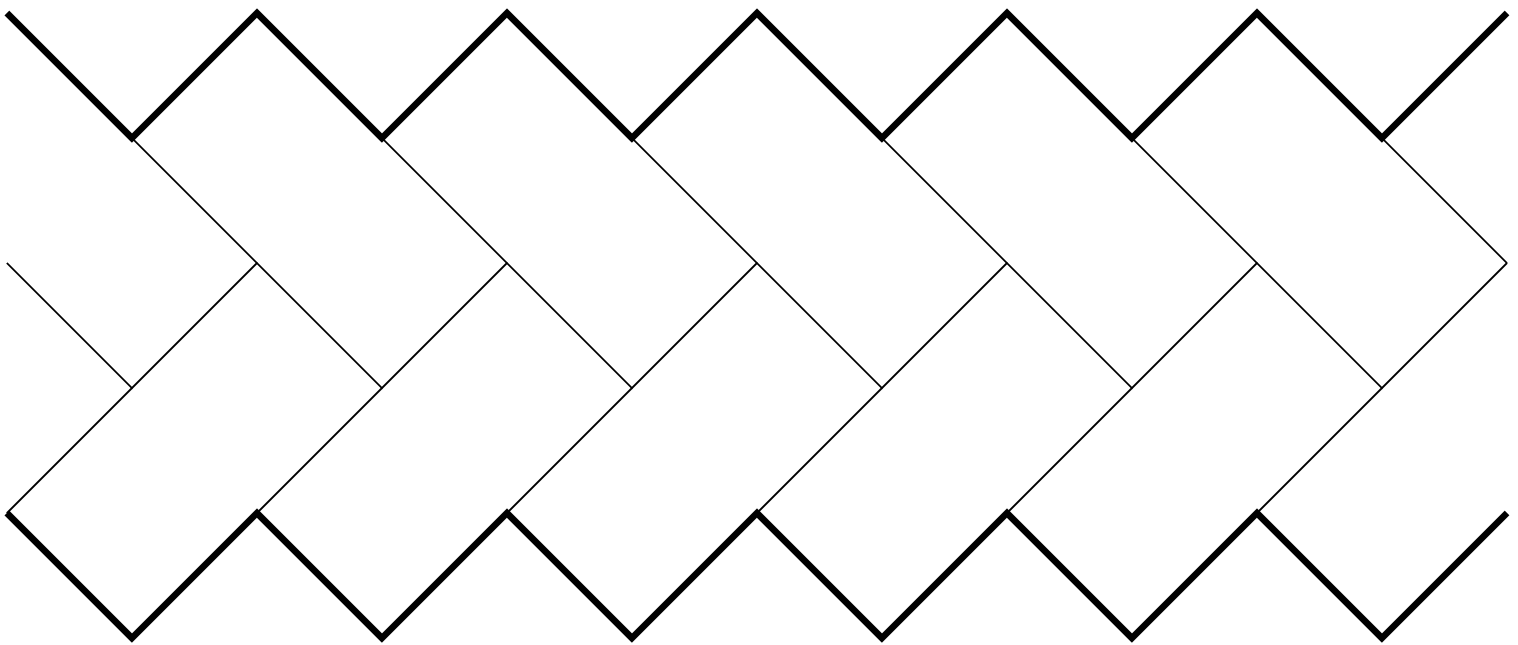}
\end{center}
\caption{Ladders}
\label{fig:ladder}
\end{figure}

Something else may happen to prevent connectivity:
the presence of {\it ladders} as in Figure \ref{fig:ladder},
where annuli have been sliced open.
Flips don't change ladders.

\begin{figure}[ht]
\begin{center}
\epsfig{height=5cm,file=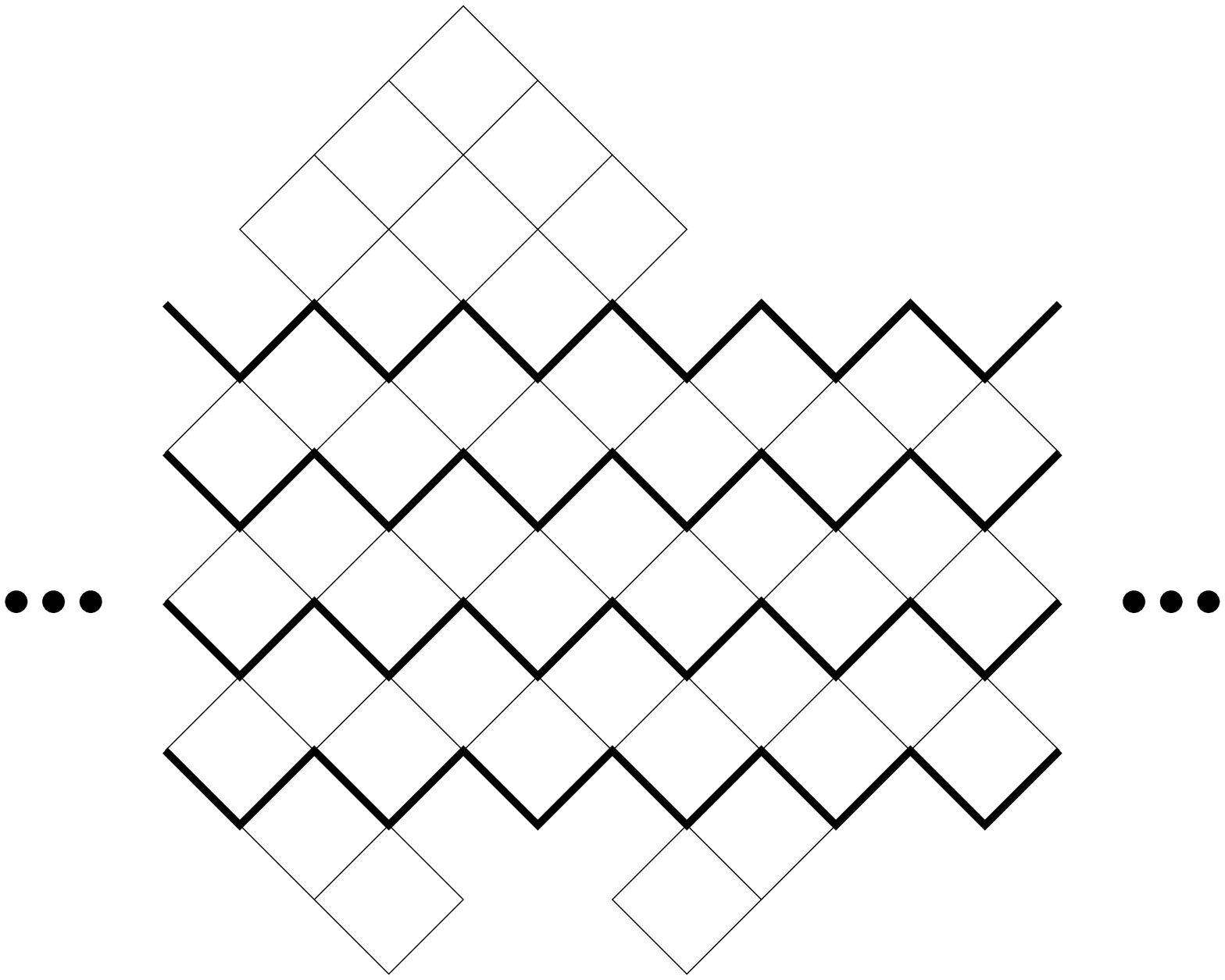} 
\end{center}
\caption{Walls}
\label{fig:zigzag}
\end{figure}

Boundaries of ladders are {\it walls}.
Tilings don't trespass walls.
Thus, for the intent of studying tilings, walls decompose annuli
into independent disks and narrow annuli (as in Figure \ref{fig:zigzag}).
There is no real loss in considering,
for the purposes of this paper, {\it wall-free annuli}.
Fortunately, there are no other obstructions to connectivity:

\begin{theo}[\cite{STCR}]
\label{theo:STCR}
Let $A$ be a wall-free quadriculated annulus.
Two tilings of $A$ can be joined by flips
if and only if they have the same flux.
\end{theo}

\section{The $q$-flux polynomial}

We will define the {\it volume} $\nu(t)$ of a tiling $t \in T_A$;
volume increases by one when applying a positive flip.
The $q$-flux polynomial $\Phi_A(p,q)$ will $(p,q)$-count tilings
with respect to flux and volume.

\begin{figure}[ht]
\begin{center}
\epsfig{height=4cm,file=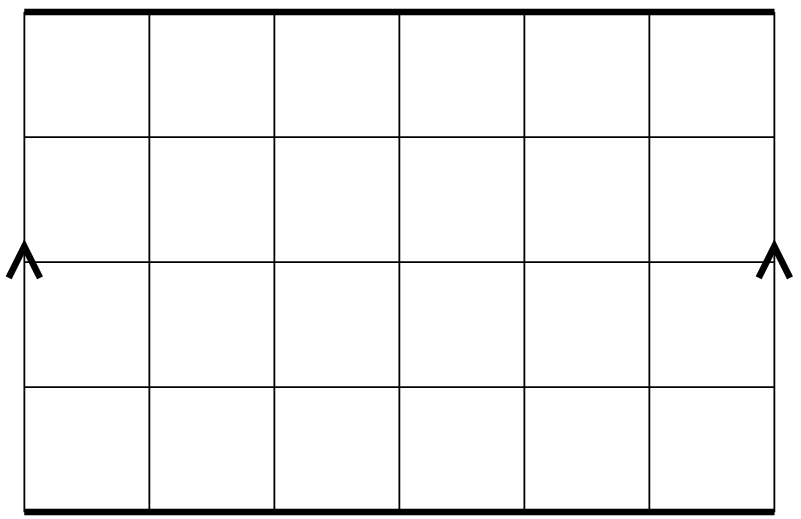}\quad
\epsfig{height=4cm,file=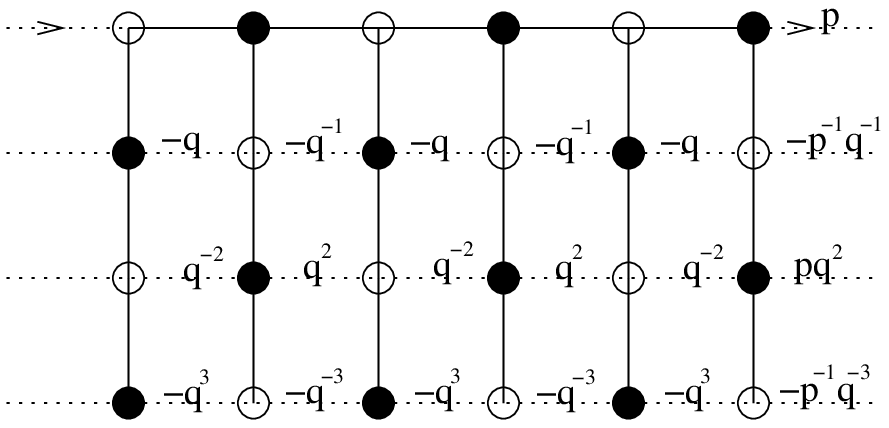}
\end{center}
\caption{Annulus and Kasteleyn weights}
\label{fig:tree}
\end{figure}

Open the annulus $A$ by a {\it cut} $\xi$ in order to obtain
a {\it track segment} $\Delta$ as in Figure \ref{fig:tree};
the annulus is obtained by identifying the vertical sides,
which both correspond to the cut $\xi$.
Draw a graph $G_A$ whose vertices are centers of squares of $A$
and whose edges join vertices of squares sharing a side.
In $G_\Delta \subset G_A$ choose a maximal tree
and assign weight 1 to its edges (solid segments in the figure).
Assign weight $p$ to an arbitrarily chosen edge in $G_A \setminus G_\Delta$
(i.e., an edge crossing the cut).
Now consider a $2 \times 2$ square contained in $A$:
say a positive flip on the square takes edges $e_0, e_1$
to edges $e_2, e_3$.
The weights $w_0, w_1, w_2, w_3$ assigned to these four edges
must satisfy $w_2 w_3 = -q w_0 w_1$:
these requirements uniquely determine weights on all edges of $G_A$.

For each tiling $t \in T_A$,
the {\it Kasteleyn weights} above obtain a signed monomial
$\pm p^{\phi(t)} q^{\nu(t)}$
by multiplying the weights of the edges asociated to dominoes in $t$.
By construction, a positive flip preserves the {\it flux} $\phi$
and increases the {\it volume} $\nu$ by one.
This also shows how to define the volume of a tiling of a quadriculated disk
(there is no flux variable).

Why are the weights signed? So that we can construct a
{\it  Kasteleyn matrix} $M_A$ whose determinant,
the {\it $q$-flux polynomial}, $(p,q)$-counts tilings, i.e.,
\[ \Phi_A(p,q) = \sum_{t \in T_A} p^{\phi(t)} q^{\nu(t)} = \pm \det(M_A). \]
Rows and colums of $M_A$ correspond to black and white squares in $A$
and entries are given by weights.
The reader should recognize the construction above as an extension
of Kasteleyn's (\cite{K}).
For the example in Figure \ref{fig:tree},
\begin{align}
\Phi_A(p,q) &= q^{-18}\,p^{-2}\,\big(\,\,
{q}^{36}\,p^4 + \notag\\
&\phantom{}
({q}^{36}+3\,{q}^{35}+3\,{q}^{34}+4\,{q}^{33}+6\,{q}^{32}+6\,{q}^{31}+7
\,{q}^{30}+6\,{q}^{29}+6\,{q}^{28}+7\,{q}^{27}+
\notag\\
&\phantom{= +}6\,{q}^{26}+6\,{q}^{25}
+7\,{q}^{24}+6\,{q}^{23}+6\,{q}^{22}+4\,{q}^{21}+3\,{q}^{20}+3\,{q}^{
19}+{q}^{18})\,p^3 + \notag\\
&\phantom{}
({q}^{30}+3\,{q}^{29}+3\,{q}^{28}+4\,{q}^{27}+9\,{q}^{26}+12\,{q}^{25}+
16\,{q}^{24}+24\,{q}^{23}+33\,{q}^{22}+
\notag\\&\phantom{= +}
41\,{q}^{21}+45\,{q}^{20}+51\,{
q}^{19}+57\,{q}^{18}+51\,{q}^{17}+45\,{q}^{16}+41\,{q}^{15}+
\notag\\&\phantom{= +}
33\,{q}^{
14}+24\,{q}^{13}+16\,{q}^{12}+12\,{q}^{11}+9\,{q}^{10}+4\,{q}^{9}+3\,{
q}^{8}+3\,{q}^{7}+{q}^{6})\,p^2 + \notag\\
&\phantom{}
({q}^{18}+3\,{q}^{17}+3\,{q}^{16}+4\,{q}^{15}+6\,{q}^{14}+6\,{q}^{13}+7
\,{q}^{12}+6\,{q}^{11}+6\,{q}^{10}+7\,{q}^{9}+
\notag\\
&\phantom{= +}6\,{q}^{8}+6\,{q}^{7}+7
\,{q}^{6}+6\,{q}^{5}+6\,{q}^{4}+4\,{q}^{3}+3\,{q}^{2}+3\,q+1)\,p + 1\,\big)
\notag
\end{align}

We are ready to state our main result.

\begin{theo}
\label{theo:ST}
Let $A$ be a balanced bicolored wall-free quadriculated annulus.
\begin{enumerate}[(a)]
\item{All non-zero roots of the polynomial $\Phi_A(p,-1)$
are roots of unity.}
\item{Let $q>0$ be fixed:
all non-zero roots of the $q$-flux
polynomial $\Phi_A(p,q)$ are distinct, negative numbers.}
\end{enumerate}
\end{theo}

In the example above,
$\Phi_A(p,-1) = p^{-2}(p^4 + p^3 + p^2 + p + 1)$ and
$\Phi_A(p,1) = p^{-2}(p^4 + 91\,p^3 + 541\,p^2 + 91\,p + 1)$
(with roots approximately equal to $-84.619$, $-6.2077$,
$-.16109$ and $-.011818$).

\section{The connection matrix}

A cut $\xi$ transforms an annulus $A$ into a track segment $\Delta$
with {\it left} and {\it right attachments}.
The $n$-th cover $A^n$ of $A$ can be obtained by juxtaposing $n$ copies
$\Delta_0, \ldots, \Delta_{n-1}$ of $\Delta$ and then {\it closing up},
i.e., identifying extreme attachments.
There are copies $\xi_{i+\half}$ of the cut between
$\Delta_i$ and $\Delta_{i+1}$,
$i = 0, \ldots, n-1$ where $\Delta_n = \Delta_0$.
A tiling of $A^n$ restricts to $\Delta_i$ yielding
a {\it tiling of the track segment} $\Delta_i$:
notice that such tilings include half-dominoes belonging to
dominoes trespassing the cuts $\xi_{i\pm\half}$
across a certain set of sides of squares.

\begin{figure}[ht]
\begin{center}
\epsfig{height=25mm,file=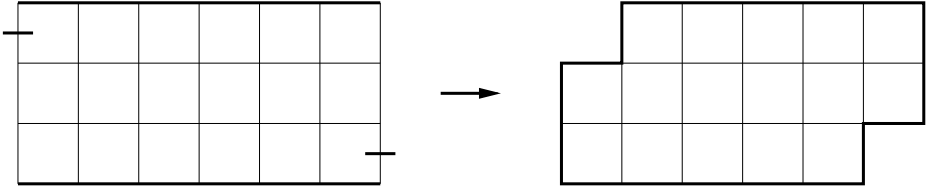}
\end{center}
\caption{Pruning a track segment to obtain a disk}
\label{fig:prune}
\end{figure}

A {\it shape} at an attachment is a set of sides contained in the attachment.
Tilings of a track segment $\Delta$ thus induce shapes on both attachments.
In particular, shape determines flux.
A pair of shapes on the attachments of $\Delta$
describes a {\it pruning} of $\Delta$
(a smaller disk or a union of disks)
by removing the squares with sides belonging to either shape,
as in Figure \ref{fig:prune}.
If two sides of the same square are selected, pruning is not defined.
Tilings of the track segment with prescribed shapes
are in natural bijection with tilings of the pruned segment.
Pairs of shapes for which pruning is not defined
are not induced by any tilings of $\Delta$.

We now construct the {\t connection matrix} $C_{\Delta}$.
Order shapes so that flux is non-decreasing.
Rows and columns of $C_{\Delta}$ correspond
to shapes at the left and right attachments.
Given a pair of shapes, the associated entry $q$-counts
tilings of the pruned segment (the entry is 0 if pruning is not defined).
We do not discuss how to compare volumes of tilings on different prunings---%
this introduces a certain ambiguity, given by possible multiplication
of entries of the matrix by different powers of $q$;
the issue is resolved via height functions in \cite{ST}.
Thus, the matrix $C_{\Delta}$ is block diagonal with blocks
$C_{\Delta,f}$ labeled by $f$, the value of the flux.
Also,
\[ \Phi_A(p,q) = p^\ast q^\ast \sum_f (\tr C_{\Delta,f}) p^f. \]
Let $\Delta^n$ be the juxtaposition of $n$ copies of $\Delta$, i.e.,
$\Delta^n$ is obtained by cutting $A^n$ across any $\xi_{i+\half}$.
The connection matrix of $C_{\Delta^n}$ is $(C_\Delta)^n$:
essentially, entries of the $n$-th power of the adjacency matrix
of a graph count paths of length $n$ between two vertices.

The polynomial $\Phi_{A^n}(p,q)$ admits two useful descriptions.
From the previous paragraph,
\[ \Phi_{A^n}(p,q) = p^\ast q^\ast \sum_f (\tr (C_{\Delta,f})^n) p^f. \]
The second relates the roots $\lambda_i$ of $\Phi_A(p,q)$
(in $p$ for fixed $q$) with the roots $\lambda_{i,n}$ of $\Phi_{A^n}(p,q)$:
\[ \lambda_{i,n} = (-1)^{n+1} \lambda_i^n. \]
This follows from constructing a Kasteleyn matrix $M_{A^n}$ from $M_A$
and diagonalizing both matrices by exploiting the
obvious action of $\ZZ/(n)$ over $A^n$.

\section{The proof}

We begin with item (a).  Suppose
that $\Delta$ is a track segment obtained from cutting $A$.
Consider, for a fixed value $f$ of the flux and $q = -1$,
the blocks $C_{\Delta^n,f} = (C_{\Delta,f})^n$
of the connection matrix $C_{\Delta^n}$ of the juxtaposition $\Delta^n$.
From the left hand side of the equality, the block
entries $-1$-count the numbers of tilings of (unions of) disks obtained
by pruning, or are equal to zero in case pruning is not defined.
From Theorem \ref{theo:DT}, these entries, for all $n$, take very few values.

Thus, for all $n$, there are only finitely many possible values for 
the matrices $C_{\Delta^n,f}$. In particular, there are powers 
$n_0$ and $n_1$ such that, for all values $f$ of the flux, we have
$C_{\Delta^{n_0},f} = C_{\Delta^{n_1},f}$, implying in turn the equality
of the polynomials $\Phi_{A^{n_0}}(p,-1) = \Phi_{A^{n_1}}(p,-1)$. 
Without loss, $n_0$ and $n_1$ can be taken to be powers of 2
so that $n_1 = n' n_0$ for a natural number $n' > 1$.
Combine this fact with the relationship 
$\lambda_{i,n} = (-1)^{n+1} \lambda_i^n$
to learn that, up to signs, raising to $n'$ induces a permutation
on the set of roots of  $\Phi_{A^{n_0}}(p,-1)$.
For a sufficiently high power $n''$ of $n'$,
raising to $n''$ keeps all roots of $\Phi_{A^{n_0}}(p,-1)$ fixed:
the nonzero roots are therefore roots of unity,
and item (a) is proved.

The proof of (b) takes a different route.
Take $q > 0$ fixed throughout the proof.
The first step is to prove that each block $C_{\Delta,f}$
has a simple eigenvalue $\Lambda_f > 0$ of absolute value larger
than that of any other eigenvalue.

A shape is {\it bi-active} if a tiling
of the bi-infinite band
$\Delta^\infty = \cdots\Delta_{-1}\Delta_0\Delta_1\cdots$
exists with the prescribed shape at the cut $\xi_\half$.
The {\it bi-active submatrix} $C_{\ast,\Delta,f}$ is the intersection of rows
and columns of $C_{\Delta,f}$ associated with bi-active shapes.
We leave it to the reader to check that the corresponding submatrix
of $(C_{\Delta,f})^n$ equals $(C_{\ast,\Delta,f})^n$
and that the spectra of $C_{\Delta,f}$ and $C_{\ast,\Delta,f}$
coincide, up to null eigenvalues;
details are given in \cite{ST}.

We now show that sufficiently large powers of $C_{\ast,\Delta,f}$
have only positive entries, i.e, that
for a given track segment $\Delta$ and a value $f$ of the flux,
there exists an integer $N$ such that for all $n > N$ and any two
bi-active shapes $\ell$ and $r$, there exists a tiling of
$\Delta^n$ with these prescribed shapes at the left and right attachments.
We proceed to join $\ell$ and $r$ by a tiling of a long track segment
$\Delta^n$.
By hypothesis, $\ell$ extends as a tiling
to the right of $\xi_\half$ in $\Delta^\infty$:
since there are only finitely many shapes,
two cuts $\xi_{I_L + \half}$ and $\xi_{I_L + P_L + \half}$ see the same shape.
By repeating this chunk of tiling, $\ell$ also extends as an eventually
periodic tiling with period $P_L$ after an initial stretch of length $I_L$.
Similarly, $r$ also extends as an eventually periodic tiling
with period $P_R$ and final stretch of length $F_R$, say,
to the left of $\xi_{n + \half}$.

Assume without loss that $P_L$ and $P_R$ are both even,
$I_L \equiv 0 \pmod {P_LP_R}$ and $F_R \equiv n \pmod {P_LP_R}$;
call $t_I$ and $t_F$ the restrictions of the infinite tilings above
to $\Delta^{I_L}$ and $\Delta^{F_R}$.
Both infinite tilings restricted to periodic stretches
(i.e., tilings of $\Delta^{P_L}$ and $\Delta^{P_R}$)
by replication give rise to tilings $t_L$ and $t_R$
of the annulus $A^{P_LP_R}$ (and thus of $\Delta^{P_LP_R}$).
It is not too hard to see that if $A$ is wall-free than so is $A^{P_LP_R}$
(see \cite{ST} for details).
From Theorem \ref{theo:STCR}, there exists a sequence of flips
connecting these tilings.
Call the tilings in this sequence $t_L = t_0, t_1, \ldots, t_M = t_R$.

We would like to join shapes $\ell$ and $r$ by a tiling obtained
by juxtaposing $t_I, t_0, t_1, \ldots, t_M, t_F$:
this is not quite correct since $t_i$ and $t_{i+1}$
may differ at the common cut. This difficulty is circumvented
by constructing tilings $t_{i+\half}$ of $\Delta^{P_LP_R}$
with the same left shape as $t_i$ and the same right shape as $t_{i+1}$
and then juxtaposing
$t_I, t_\half, t_{1.5}, \ldots, t_{M - \half}, t_F$.
More precisely, if $t_i$ and $t_{i+1}$ have the same left shape
take $t_{i+\half} = t_i$;
otherwise, take $t_{i+\half}$ to coincide with $t_i$ on the left
$\Delta^{P_LP_R/2}$ subsegment and with $t_{i+1}$ on the right subsegment.
This obtains a tiling in $\Delta^{I_L + (M - 1)P_LP_R + F_R}$:
since this number is congruent to $n$ modulo $P_LP_R$,
the desired tiling can be obtained, for sufficiently large $n$,
by inserting copies of $t_M$ before $t_F$.

The Perron-Frobenius theorem (\cite{G}) applied to $C_{\ast,\Delta,f}$
then completes the proof of the first step.
 
Label the nonzero roots of $\Phi_A(p,q)$ as
$|\lambda_1| \ge |\lambda_2| \ge \ldots \ge |\lambda_m| > 0$.
Assume by induction that $\lambda_1, \ldots, \lambda_{k-1}$
are real negative and that
\[ |\lambda_1| > \cdots > |\lambda_{k-1}| > |\lambda_k| = \cdots
= |\lambda_{k'}| > |\lambda_{k'+1}|, \]
with $k' \ge k$.
We must prove that $k = k'$ and that $\lambda_k$ is real negative.
Consider the usual symmetric function
$$\sigma_k(x_1,\ldots,x_m) =
\sum_{1 \le i_1 < \cdots < i_k \le m} x_{i_1} \cdots x_{i_k}.$$
From what we have seen above, 
\[ \tr C_{\Delta,f_{\max} - k} =
a^n \sigma_k((-\lambda_1)^n, \ldots, (-\lambda_m)^n) =
(1 + o(1)) \Lambda_{f_{\max} - k}^n \]
when $n$ goes to infinity.
Here $a\,p^{f_{\max}}$ is the leading monomial of $\Phi_A(p,q)$.

The expression $ \sigma_k((-\lambda_1)^n, \ldots, (-\lambda_m)^n) $
is the sum of $k' - k + 1$ terms of the form
$((-1)^k \lambda_1\lambda_2\cdots\lambda_{k-1}\lambda_\ell)^n$,
$\ell = k, \ldots, k'$ and other terms
which grow at exponentially smaller rates.
Thus
\[ (-1)^{kn} \lambda_1^n\lambda_2^n\cdots\lambda_{k-1}^n
(\lambda_k^n + \cdots + \lambda_{k'}^n) = (1 + o(1)) (\Lambda_k/a)^n \]
whence
\[ b_k^n + \cdots + b_{k'}^n = 1 + o(1) \]
where $b_\ell = - \lambda_\ell/|\lambda_\ell|$.
The upshot is that $b_k, \ldots, b_{k'}$ belong to the unit circle
and we can take arbitrarily large $n$ such that 
$2 (k' - k + 1)|b_\ell^n - 1| < 1$ for all $\ell$
and therefore
\[ |(k' - k + 1) - ( b_k^n + \cdots + b_{k'}^n )| < 1/2, \]
a contradiction unless $k = k'$.
Finally, $b_k = 1$ and $\lambda_k$ is real negative.
This concludes the proof of item (b).

As an application of (b), we state the result below.
A sequence $a_k$ of non-negative real numbers is {\it log-concave} if
$a_k^2 \ge a_{k-1} a_{k+1}$ for all $k$.
In particular, log-concave sequences are either monotone or unimodal.

\begin{coro}
Let $q$ be a fixed positive real number and
let $a_f$ be the coefficient of $p^f$ in $\Phi_A(p,q)$.
Then the sequence $a_f$ is log-concave.
\end{coro}


{\obeylines
Nicolau C. Saldanha and Carlos Tomei
Depto. de Matem{\'a}tica, PUC-Rio
R. Mq. de S. Vicente 225
Rio de Janeiro, RJ 22453-900, Brazil
nicolau@mat.puc-rio.br
tomei@mat.puc-rio.br
}

\end{document}